\documentclass[a4paper,11pt,reqno]{amsart}
\pdfoutput=1

\usepackage[]{amsmath}
\usepackage{amsfonts}
\usepackage{amssymb}
\usepackage{xcolor}
\usepackage{mathrsfs}

\usepackage{tikz}

\usepackage[margin=3cm]{geometry} 
\usepackage{hyperref}
\catcode`\@11\relax
\newif\ify@autoscale \y@autoscaletrue \def\Yautoscale#1{\ifnum #1=0
  \y@autoscalefalse\else\y@autoscaletrue\fi}
\newdimen\y@b@xdim
\newdimen\y@boxdim \y@boxdim=13pt
\def\Yboxdim#1{\y@autoscalefalse\y@boxdim=#1}
\newdimen\y@linethick    \y@linethick=.3pt
\def\Ylinethick#1{\y@linethick=#1}
\newskip\y@interspace \y@interspace=0ex plus 0.3ex
\def\Yinterspace#1{\y@interspace=#1}
\newif\ify@vcenter   \y@vcenterfalse
\def\Yvcentermath#1{\ifnum #1=0 \y@vcenterfalse\else\y@vcentertrue\fi}
\newif\ify@stdtext   \y@stdtextfalse
\def\Ystdtext#1{\ifnum #1=0 \y@stdtextfalse\else\y@stdtexttrue\fi}
\newif\ify@enable@skew   \y@enable@skewfalse
\expandafter\ifx\csname enableskew\endcsname\relax
 \y@enable@skewfalse \else \y@enable@skewtrue\fi
\def\y@vr{\vrule height0.8\y@b@xdim width\y@linethick depth 0.2\y@b@xdim}
\def\y@emptybox{\y@vr\hbox to \y@b@xdim{\hfil}}
\ify@enable@skew
 \def\y@abcbox#1{\if :#1\else
   \y@vr\hbox to \y@b@xdim{\hfil#1\hfil}\fi}
 \def\y@mathabcbox#1{\if :#1\else
   \y@vr\hbox to \y@b@xdim{\hfil$#1$\hfil}\fi}
\else
 \def\y@abcbox#1{\y@vr\hbox to \y@b@xdim{\hfil#1\hfil}}
 \def\y@mathabcbox#1{\y@vr\hbox to \y@b@xdim{\hfil$#1$\hfil}}
\fi
\def\y@setdim{%
  \ify@autoscale%
   \ifvoid1\else\typeout{Package youngtab: box1 not free! Expect an
     error!}\fi%
   \setbox1=\hbox{A}\y@b@xdim=1.6\ht1 \setbox1=\hbox{}\box1%
  \else\y@b@xdim=\y@boxdim \advance\y@b@xdim by -2\y@linethick
  \fi}
\newcount\y@counter
\newif\ify@islastarg
\def\y@lastargtest#1,#2 {\if\space #2 \y@islastargtrue
  \else\y@islastargfalse\fi}
\def\y@emptyboxes#1{\y@counter=#1\loop\ifnum\y@counter>0
  \advance\y@counter by -1 \y@emptybox\repeat}
\def\y@nelineemptyboxes#1{%
  \vbox{%
    \hrule height\y@linethick%
    \hbox{\y@emptyboxes{#1}\y@vr}
    \hrule height\y@linethick}\vskip-\y@linethick}
\def\yng(#1){%
  \y@setdim%
  \hskip\y@interspace%
  \ifmmode\ify@vcenter\vcenter\fi\fi{%
  \y@lastargtest#1,
  \vbox{\offinterlineskip
    \ify@islastarg
     \y@nelineemptyboxes{#1}
    \else
     \y@ungempty(#1)
    \fi}}\hskip\y@interspace}
\def\y@ungempty(#1,#2){%
  \y@nelineemptyboxes{#1}
  \y@lastargtest#2,
  \ify@islastarg
   \y@nelineemptyboxes{#2}
  \else
   \y@ungempty(#2)
  \fi}
\def\y@nelettertest#1#2. {\if\space #2 \y@islastargtrue
  \else\y@islastargfalse\fi}
\def\y@abcboxes#1#2.{%
  \ify@stdtext\y@abcbox#1\else\y@mathabcbox#1\fi%
  \y@nelettertest #2.
  \ify@islastarg\unskip%
   \ify@stdtext\y@abcbox{#2}\else\y@mathabcbox{#2}\fi%
  \else\y@abcboxes#2.\fi}
 \newdimen\y@full@b@xdim
 \newcount\y@m@veright@cnt
\ify@enable@skew
 \def\y@get@m@veright@cnt#1#2.{%
   \if :#1 \advance\y@m@veright@cnt by 1\y@get@m@veright@cnt#2.\fi}
 \let\y@setdim@=\y@setdim
 \def\y@setdim{%
   \y@setdim@ \y@full@b@xdim=\y@b@xdim
   \advance\y@full@b@xdim by 1\y@linethick}
 \def\y@m@veright@ifskew#1{
   \y@m@veright@cnt=0 \y@get@m@veright@cnt#1.
   \moveright \y@m@veright@cnt\y@full@b@xdim}
\else
 \def\y@m@veright@ifskew#1{}
\fi
\def\y@nelineabcboxes#1{%
  \y@nelettertest #1.
  \ify@islastarg
   \y@m@veright@ifskew{#1}
    \vbox{
      \hrule height\y@linethick%
      \hbox{\ify@stdtext\y@abcbox#1\else\y@mathabcbox#1\fi\y@vr}
      \hrule height\y@linethick}\vskip-\y@linethick
  \else
   \y@m@veright@ifskew{#1}
    \vbox{
      \hrule height\y@linethick%
      \hbox{\y@abcboxes #1.\y@vr}%
      \hrule height\y@linethick}\vskip-\y@linethick
  \fi}
\def\young(#1){%
  \y@setdim%
  \hskip\y@interspace%
  \y@lastargtest#1,
  \ifmmode\ify@vcenter\vcenter\fi\fi{%
  \vbox{\offinterlineskip
    \ify@islastarg\y@nelineabcboxes{#1}%
    \else\y@ungabc(#1)%
    \fi}}\hskip\y@interspace}
\def\y@ungabc(#1,#2){%
  \y@nelineabcboxes{#1}%
  \y@lastargtest#2,
  \ify@islastarg\y@nelineabcboxes{#2}%
  \else\y@ungabc(#2)%
  \fi}
\catcode`\@12\relax
 
\newtheorem{thm}{Theorem}[section]

\newtheorem{lem}[thm]{Lemma}
\newtheorem{cor}[thm]{Corollary}

\theoremstyle{definition}

\def\Qbar{ \overline{Q} }

\def\obar{ \overline{1} }
\def\zbar{ \overline{0} }

\def\om{ \overline{m} }

\def\one{\hbox{{1}\kern-.25em\hbox{l}}}
\def\semiplus{\raisebox{-.12pt}{\scalebox{1.2}{$\supset$}}\hskip-13pt}
\def\mh{|\mbox{\hskip.2cm}\rangle}
\def\nn{\nonumber}

\definecolor{vertclair}{RGB}{0,153 0}

\definecolor{gris}{RGB}{153,153,153}

\title[Indecomposable doubling] 
 {Indecomposable doubling for representations of the type I  Lie superalgebras
$sl(m/n)$ and $osp(2/2n)$.}
\author[Jarvis and Thierry-Mieg]{Peter D. Jarvis${}^1$ and Jean Thierry-Mieg${}^2$}
\address{${}^1$School of Natural Sciences (Mathematics and Physics), University of Tasmania, PO Bag 37, Hobart, Tasmania 7001, Australia
(Alexander von Humboldt Fellow); 
\mbox{${}^2$NCBI}, National Library of Medicine, National Institute of Health, 8600 Rockville Pike, Bethesda MD20894, U.S.A.
}
\email{${}^1$peter.jarvis@utas.edu.au, ${}^2$mieg@ncbi.nlm.nih.gov}

\date{\today }
\begin{document}

\begin{abstract}
We establish that for the type I Lie superalgebras $sl(m/n)$ and $osp(2/2n)$, each Kac module admits a 1 parameter family of indecomposable double extensions. The result follows from the explicit evaluation of the $H^1$ Lie superalgebra cohomology valued in the tensor product of the module and its dual. 
\end{abstract}
\maketitle
\vfill
\vfill

\pagebreak

\section{Introduction} 
\mbox{}

\noindent
With its motivation in physics and origins in Bose-Fermi statistics, and applications in quantum mechanics and condensed matter as well as supersymmetry in particle physics, the analysis of Lie superalgebras was formalized in a mathematical framework with 
Kac's \cite{Kac1977classification} classification of finite-dimensional simple Lie superalgebras, and 
subsequently with the identification of finite dimensional irreducible representations and certain character formulae \cite{Kac1977characters} (see also Scheunert \cite{scheunert2006theory}).
For a review of fundamentals of superalgebra theory in the context of applications, see for example \cite{frappat2000dictionary}\, or the text \cite{cornwell1992group}\,. 

However, even in finite dimensions, and despite their unavoidable involvement in the resolution of tensor products of irreducible representations, the analysis of indecomposable representations of Lie superalgebras remains in general a difficult open problem.
Initial studies focussed on low rank cases, such as 
$osp(1/2)\cong B(0,1)$\, in Kac's classification, as well as $sl(2/1)\cong A(1,0)$, or special cases such as $osp(1/2n) \cong B(0,n)$\, (notational details will be given in the main text below). For example, in \cite{ScheunertNahmRittenberg1977irreps,marcu1980tensor}\,,
the indecomposability of the tensor product 
of two quartet representations of $sl(2/1)$ was noted, and in \cite{marcu1980tensor, marcu1980representations,Su1992,gotz2007representation}\,, broader studies were undertaken on the identification of arbitrary (finite dimensional) representations, as well as the identification of indecomposables in the resolution of tensor products. For an analysis of indecomposable representations of $sl(n/1)$ see \cite{germoni1997representations}\, and for other cases 
\cite{germoni1998indecomposable,germoni2000indecomposable}. 
For the role of indecomposables in the construction of so-called atypical irreducibles, and their characters, we cite \cite{Leites1975cohomology,Serganova1996KL,SuZhang2007character}\,.

The present paper presents a result on a specific form of indecomposable representation for a class of the so-called type I simple Lie superalgebras, 
namely the special linear $sl(m/n) \cong A(m-1,n-1)$\, ($m\ne n \ge 1$\,), and orthosymplectic superalgebras $osp(2/2n) \cong C(n)$\,.
We establish that any Kac module (including indecomposable cases thereof) can be extended indecomposably to a one parameter family of doubled modules, of twice the dimension. 

This result is, in its own right, a striking property of superalgebra representations. But it is potentially of  importance in physical applications, wherein states of a given system
occur as repeated copies (or, possibly, as multiple copies). It has been proposed for a number of years (\cite{ne1979irreducible,fairlie1979higgs,dondi1979supersymmetric,neeman1980geometrical,thierry1982exterior}\,; see also \cite{ThierryMieg2021newchi,ThierryMiegJarvis2021newchi}\,), that the $SU(2/1)$ superalgebra could play a role in the standard model of the fundamental interactions.
The lepton and quark states, graded by their chiralities, are precisely described by quartet representations of $SU(2/1)$.
Furthermore, generational repetition and mixing of the families of quarks and leptons, corresponding to the electron and the muon, can be accommodated by the existence
\cite{marcu1980representations,Su1992} of 
indecomposable doubles of the corresponding $SU(2/1)$ quartet. Thus these indecomposable structures offer
an algebraic understanding of the weak interaction Cabibbo mixing angle,
with extensions to the Cabibbo-Kobayashi-Maskawa three-generation mixing (
\cite{coquereaux1991elementary,ScheckHaussling1998triangular,haussling1998leptonic}\,; see also \cite{ThierryMiegJarvis2021indec})\,.

In section \ref{sec:KacModules} below, we give a brief technical background to the identification and 
construction of the Kac modules for $sl(m/n)$ and $osp(2/2n)$.
In section \ref{sec:Relative} the method used for identification of module extensions, 
namely the analysis of appropriate (co)homology groups,
is reviewed, in the context of the salient
 literature. A useful lemma on the construction of induced modules
is proven (Lemma \ref{lem:DoubleSided} below), and Appendix \ref{sec:DoubleSided}\, provides further details of the induced module construction 
in relation to this. This is then used (section \ref{sec:Relative}) to produce a count and 
analysis of chain operators, from which the results follow 
(Theorem \ref{thm:IndecDoubling}, Corollary \ref{cor:UniqueIndecDoubling}). In the concluding section \ref{sec:Discussion} additional discussion and some generalizations are given,
via examples arising from the simplest $sl(2/1)$ case. For this specific
case, further details are also given, in relation to its application to mixing in the standard model
formulated via $SU(3)\times SU(2/1)$\,\cite{ThierryMiegJarvis2021newchi,ThierryMiegJarvis2021indec}\,.

\section{Kac modules for type I LSA and ``double'' induction}
\label{sec:KacModules}
\mbox{}
While irreducible representations of Lie superalgebras can be constructed using a graded version of the standard
Verma module induction via suitable Borel superalgebras based on the Cartan subalgebra, we focus here on the construction
of Kac \cite{Kac1977classification}, whereby the induction step builds on (parabolic) sub superalgebras involving the even subalgebra. Here we briefly introduce the notation and outline the constructions, with emphasis on structural relationships
between modules so constructed and their duals.

A Lie superalgebra  is a (finite dimensional) ${\mathbb Z}_2$-graded linear space $L= L_{\zbar} + L_{\obar}$ comprising even and odd subspaces wherein $|[X,Y]| = |X| + |Y|\,\mbox{mod 2}$ and for the graded Lie bracket
\begin{align}
\label{eq:LSAdef}
[X, Y] =&\, -(-1)^{|X||Y|}[Y,X]\,, \nn \\ 
[X,[Y,Z]] + (-1)^{|X|(|Y|+|Z|)}[Y, [Z,X]]&\, + (-1)^{|Z|(|X|+|Y|)}[Z,[X,Y]] = 0\,,
\end{align}
with $X,Y,Z$\, homogeneous.
For the type I Lie superalgebras of concern here, we have
$L_{\zbar} = L_0 +{\mathbb C}$ and $L_{\obar} = L_{-1} + L_{+1}$ with $L_0$  semisimple and 
${\mathbb C} \cong {\mathbb C} Y$ a one-dimensional
abelian Lie algebra, spanned by a generator $Y$\, and where in the direct sum we have
${\{}L_{\pm 1}, L_{\pm 1}{\}} \equiv 0$\,, that is 
\begin{align}
 [Y, \Qbar] =&\, +\Qbar\,,\nn \\
[Y,Q] = &\, -Q\,,\nn \\
\mbox{and}\qquad {\{}Q, Q{\}} =&\, 0 = {\{}\Qbar, \Qbar{\}}\,, 
\end{align}
where $\Qbar \in L_{+1}$  and $Q \in L_{-1}$\,, respectively, with ${\{}\Qbar, Q{\} }\in L_{\zbar}$\,. (Here for homogeneous elements, ${\{} \,\cdot\,, \,\cdot\,{\}}$ is the symmetric (anticommutator)
symbol for brackets of odd generators\,; brackets with even elements or even-odd elements involve antisymmetric (commutator) symbols ${[} \,\cdot\,, \,\cdot\,{]}$\,.)  We now turn to the above-mentioned Kac module construction.

In the type I case, the superalgebra $L$ admits natural induced constructions via the semidirect product sub-superalgebras $L_{\pm}= L_{\zbar}+ L_{\pm 1}$\,. Taking an even $L_{\zbar}$-module $U$, one defines
a corresponding  $L_{}$-module via the `highest weight' condition
\[
L_{+1} U =0\,, \quad \mbox{viz.} \quad \Qbar u =0\,, \quad \forall\,\, \Qbar \in L_{+1} \,\,, u \in U\,,
\]
and similarly for $L_{-1}$\, wherein a lowest weight condition is applied. The respective induced modules are
\begin{align}
\overline{V_+(U)}{} = &\, \mbox{Ind}\! \left.\right|_{L_{+}}^L (U)\, := {\sf U}(L) \otimes_{L_{+}} U\,, \nn \\
\overline{V_-(U)}{} = &\, \mbox{Ind}\! \left.\right|_{L_{-}}^L (U)\,:= {\sf U}(L) \otimes_{L_{-}} U\,.
\end{align}
This construction yields irreducible or decomposable $L$-modules, depending on whether certain conditions are satisfied by the highest weight vector of $U$. 
In the former case, the modules are called typical. Evidently, in view of the Poincar\'{e}-Birkhoff-Witt theorem, the restricted tensor product construction means that a spanning set is derived solely from tensoring by
the enveloping algebra of $L\backslash L_{\pm}$\,, namely the exterior algebra
\smash{$\scalebox{1.2}{\boldmath{$\wedge$}}(L_{\mp })$}\,, respectively, so that typical irreducible modules have dimension $2^{\frac 12|L_{\obar}|} \mbox{dim}(U)$\,.  Otherwise the Kac module (of the same dimension) is indecomposable, and further atypical irreducibles are recovered by factoring by the maximal invariant submodule (and denoted $V_\pm(U)$\,). 
For the present discussion we do not need further structural details from the classification of simple Lie superalgebras or their representations; specific cases and examples are treated in section \ref{sec:Discussion} below. We now turn to two useful lemmas regarding dual modules, which will facilitate the identification of the indecomposable doubled modules in the next section.

\begin{lem}
[Duals as lowest weight modules]\mbox{}\\
The dual of a Kac module $\overline{V_+(U)}{}$ is isomorphic to the opposite induced module
$\overline{V_-(U{}^*)}{}$: that is,
\begin{align}
\Big(\mbox{Ind}\! \left.\right|_{L_{+}}^L (U)\Big)^* \cong 
\mbox{Ind}\! \left.\right|_{L_{-}}^L (U^*)\,.
\end{align}
\label{lem:DualKac}\end{lem}
\begin{lem}
[``Double'' induction]\mbox{}\\
The tensor product of a Kac module $\overline{V_+(U)}$ and its dual $\overline{V_-(U{}^*)}$ is isomorphic to the $L$-module induced from the $L_{\zbar}$-module $U \otimes U^*$\,; that is, 
\begin{align}
\overline{V_+(U)}\otimes \overline{V_-(U{}^*)} \cong \mbox{Ind}\! \left.\right|_{L_{\zbar}}^L \big(U \otimes U^*\big)\,.
\end{align}
\label{lem:DoubleSided}\end{lem}
\noindent
\textbf{Proofs:}\\
Evidently, $\overline{V_+(U)}{}^*$ is completely reducible to a sum of $L_{\zbar}$-modules which are the duals of the corresponding components of $\overline{V_+(U)}$ with opposite weight. If $v$ is an element of such an $L_{\zbar}$-module with weight $\lambda$, then for duals $w^*$, we have 
\[
h^* \cdot w^* (v) = w^*(-h\cdot v) = -\lambda(v)w^*(v)\,,
\]
so that $h^*\cdot w^* = -\lambda(v)w^*$ for $w^*(v)\ne 0$\,. If $U$ has highest weight $\Lambda$\,, then $U^*$ has lowest weight $-\Lambda$\,, and given
$\Qbar\cdot u =0$ for $u\in U$, we have $\Qbar{}^* \cdot u^* =0$ and indeed the coadjoint 
generators have the opposite ${\mathbb Z}$-grading by virtue of 
$[Y^*, \Qbar{}^*]=  -\Qbar{}^*\,$, and $[Y^*, Q{}^*]=  +Q{}^*\,$\,. In this way we realize 
\begin{eqnarray*}
\overline{V_+(U)}{}^*&\cong& \overline{V_-(U{}^*)} 
\end{eqnarray*}
\mbox{}\\[-.8cm]\mbox{}\hfill $\Box$

\noindent
The second lemma is established via the mapping 
\begin{align}
\label{eq:LinSpaceIso}
\pi:\big(m \otimes_+ u \big)\otimes 
\big(\overline{m}{} \otimes_- u^* \big)
\mapsto &\, \big(m \overline{m}{}\big) \otimes \big( u \otimes u^* \big)
\end{align}
where we have written $\otimes_\pm$ for $\otimes_{L_\pm}$\,).
Here $m , \, \overline{m}{} $
are antisymmetrized monomials in $Q$\, $\Qbar$\,, such that
via the Poincar\'{e}-Birkhoff-Witt theorem, the ordered elements $m \overline{m}{}$ provide a spanning set for the ``double'' induced module\footnote{Taking the grading of the module $U$ to be even.} so that $\pi$ is a linear space isomorphism. The equivalence of superalgebra modules can be established by a careful 
analysis of the respective actions (details are given appendix \ref{sec:DoubleSided})\,.\mbox{}\hfill $\Box$\\

\noindent
We now turn to the implications of these structural features of Kac modules for 
establishing the existence of the desired indecomposable doubles.  Concrete examples are discussed in the concluding section \ref{sec:Discussion} below.

\section{(Co)homology: indecomposable doubling for Kac modules}
\label{sec:Relative}
The classification and construction of arbitrary (finite dimensional) indecomposable representations for arbitrary Lie superalgebras is, in general, combinatorially very challenging, and remains an open problem. Studies of $sl(2/1)$ can be found in \cite{marcu1980tensor,marcu1980representations,Su1992,gotz2007representation}\,; see 
\cite{germoni1997representations,germoni1998indecomposable,germoni2000indecomposable}\,  for extensions to $sl(n/1)$\, and other cases.  Here we take advantage of the special structure of the Kac module realization outlined in the previous section,
in incorporating it into basic results relating to module extensions and 
(co)homology.

The general setup for indecomposably combining two $L$-modules $U$ and $V$\, in the superalgebra case \cite{BenAmorPinczon1991}
entails the direct sum module $W = U+V$,  the extension of $U$ by $V$, in such a way that there is an exact sequence
\[
 0 \rightarrow U \stackrel{\mu}{\scalebox{1.2}{$\rightarrow$}} W \stackrel{\varepsilon}{\scalebox{1.2}{$\rightarrow$}} V\rightarrow 0 
\]
where representations $L\rightarrow \mbox{Hom}\big(U\big)$\,, $L\rightarrow  \mbox{Hom}\big(V\big)$ are augmented by representations mapping $L\rightarrow  \mbox{Hom}\big(U,V\big)$, that is
\[
X \cdot T = X_V \cdot T -(-1)^{|X||T|}  T \cdot X_U\,.
\]
Thus, the space $U$ is embedded in $W$ as an invariant subspace, and the factor space $W/U$ is isomorphic to the module $V$\,. For the purposes of discussion of examples in section \ref{sec:Discussion} below, such an extension  will be denoted ``$V \semiplus +\hskip3pt U$'' \,.

The abstract framework for identifying the $\mbox{Hom}\big(U,V\big)$ representations which enumerate independent extensions is that of superalgebra (co)homology,
a generalization of standard Lie algebra (co)homology (see for example \cite{ChevalleyEilenberg1948}) to the Lie superalgebra case
\cite{Leites1975cohomology,ScheunertZhang1998,BenAmorPinczon1991}. Here we provide key results
\cite{BenAmorPinczon1991,tanaka1995homology} enabling the solution to the doubling problem to be reduced to concrete
evaluations. Rather than give a full presentation, we focus
on the specifics relevant to the present case. Detailed
explanations of the notation and elaboration of the mappings to be considered 
are deferred until after the statement of the main results, as follows:

\begin{thm}[Module Extensions \cite{BenAmorPinczon1991}]\mbox{}\\
The indecomposable extensions of  $L$-module $U$ by $L$-module $V$ are in correspondence with elements of the first cohomology group of $L$ valued in $Hom(V,U)$\,, namely $H^1\big(L, \mbox{Hom}(V,U)\big)$\,.
\end{thm}
\begin{lem}
[Shapiro's lemma \cite{tanaka1995homology}]\mbox{}\\
The homology of an $L$-module $V$ of a Lie superalgebra $L$ which is induced from a sub-superalgebra $L'$\,, and $L'$-module $U$, is isomorphic
to the homology of $L'$ valued in $U$:
\begin{align}
H_n(L,V) \cong&\, H_n(L',U)\,; \quad \mbox{\protect {\textrm{that is,}}} \nn \\
H_n\big(L,\mbox{\protect {\textrm{Ind}}}\! \left.\right|_{L'}^L \big(U\big)\,\big) \cong &\, H_n(L',U)\,,
\end{align}
provided the $L$-modules occurring in the associated chains are completely reducible as $L'$-modules.
\label{lem:Shapiro}
\end{lem}

\noindent
Here we state results for homology rather than cohomology groups, which are canonically dual, and in fact isomorphic for self-dual modules, as in the present case
(the distinction will arise in section \ref{sec:Discussion} below).
Given that $\mbox{Hom}(V,V) \cong V \otimes V^*$\,, Lemmata \ref{lem:Shapiro}\,, \ref{lem:DualKac} and \ref{lem:DoubleSided}  lead to a powerful simplification to the
problem of enumerating doubles of the Kac module $\overline{V_+(U)}$. We have
\begin{thm}[Indecomposable doubling]\mbox{}
\begin{align}
H_1\big(L, \overline{V_+(U)}\otimes \overline{V_-(U{}^*)}\big) \cong H_1\big(L_{\zbar}, U\otimes U^*\big)^{(L_{\zbar})}\,.
\end{align}
\mbox{}\hfill $\Box$\\[-1.5cm]
\label{thm:IndecDoubling}
\end{thm}
\mbox{}\noindent
\begin{cor}[Unique Indecomposable doubling]\mbox{}\\
We have 
\begin{align}
\mbox{dim}\big( H_1\big(\overline{V_+(U)}\otimes \overline{V_-(U{}^*)}\big) \big) =1\,. 
\end{align}
That is, each Kac module $\overline{V_+(U)}$ of the class I simple Lie superalgebras 
$sl(m/n)$ and $osp(2/2n)$ admits a 1-parameter family of indecomposable self-extensions\,. 
\mbox{}\hfill $\Box$
\label{cor:UniqueIndecDoubling}
\end{cor}

\noindent
With Theorem  \ref{thm:IndecDoubling}\,, the analysis is reduced to a straightforward examination of cases based on the identification of even subalgebra invariants arising in the relevant chain spaces (the (co)homology is $L_{\zbar}$ invariant). The details of the 
formalism needed for our analysis are as follows.
As $L_{\zbar}$ modules, the chain spaces are:
\begin{align}
B_2:= &\, \big(L_{\zbar} \wedge L_{\zbar} \big) \otimes (U \otimes U^*)\,; \nn \\
B_1 := &\, L_{\zbar} \otimes (U \otimes U^*)\,; \nn \\
B_0 := &\,  (U \otimes U^*)\,,
\end{align}
with chain maps
\begin{align}
\partial_1: B_2 \rightarrow B_1\,,\quad \partial_{0}: B_1 \rightarrow B_0\,. \nn
\end{align}
We require the first ($n=1$) homology group,
\begin{align}
H_1\big(U \otimes U^*\big) = 
\mbox{Ker}\big(\partial_{0}\big)/\mbox{Im}\big(\partial_1\big)\,.
\end{align}
The chain maps are defined in the superalgebra case, 
\begin{align}
\label{eq:BoundaryMapDef}
\partial_1\big(X\wedge Y \otimes v\big) := &\, X \otimes Yv - (-1)^{|X||Y|}Y \otimes Xv -{[}X,Y{]} \otimes v\,;\nonumber \\
\partial_0 \big(X\otimes v\big) := &\, Xv\,,
\end{align}
for homogeneous elements, with analogous definitions for the cochains in the case of the cohomology \cite{ScheunertZhang1998,tanaka1995homology}\,.

As noted, the (co)chain maps intertwine the Lie superalgebra action. Here we require only the analysis of the even $L_{\zbar}$ homology, and
by Theorem \ref{thm:IndecDoubling}\,, 
\begin{align}
H_1\big(U \otimes U^*\big) \cong 
\mbox{Ker}\big(\partial_{0}{}^{(L_{\zbar})}\big)/\mbox{Im}\big(\partial_1{}^{(L_{\zbar})}\big)
\end{align}
where $\partial_n{}^{(L_{\zbar})}$\,, $n=0,1$ are the restriction of the chain maps to the $L_{\zbar}$ invariants $B_n{}^{(L_{\zbar})}$\,, $n=0,1$\,. 
\\

\noindent
\textbf{Proof of Corollary \ref{cor:UniqueIndecDoubling}}\\
The $L_{\zbar}$ module $(U \otimes U^*)$ is completely reducible to a direct sum
\begin{align}
(U \otimes U^*)\cong W_0 \oplus W_{ad} \oplus \sum \oplus W
\end{align}
of $L_{\zbar}$- modules, where $W_0:= \langle w_0\rangle$ is (by Frobenius' theorem) the unique singlet, $W_{ad}$ is the (isotypical) $L_{\zbar}$-adjoint subspace, and $\sum \oplus W$ represents the remaining decomposition into
non-adjoint $L_{\zbar}$ irreducible modules. We note that
for the semisimple Lie algebras $L_0$\,, $sl(n)$ and $sp(2n)$\,, the antisymmetric square of the adjoint representation contains the adjoint with multiplicity 1 so that the isotypical space $W_{ad}$ is in fact irreducible. We proceed by introducing an explicit basis
$\{ J^a \}$ of generators for the Lie algebra for $L_0$, with commutation relations
$[ J^a, J^b] = f^{ab}{}_c J^c$\,, and corresponding basis $\{ w^a \} $ for the adjoint
module $W_{ad}$\,, a constituent of the reduction of $(U \otimes U^*)$\,, with adjoint action $J^a w^b = -f^{ab}{}_c w^c$\,. Here we take, for the semisimple Lie algebra $L_0$\,, antisymmetric structure constants with Killing form proportional to $\delta_{ab}$\,. Finally let $\{ w^i \}$ be  a basis for one of the  non-adjoint irreducible modules $W \subset (U \otimes U^*)$\,. We enumerate the $L_{\zbar}$ invariants spanning 
the \smash{$B_n^{(L_{\zbar})}$} chain spaces, $n=0,1,2$\,, as 
\begin{align}
I_0 =&\, w_0\,;\nn \\
I_1 =&\, J^a \otimes w^a\,; \nn \\
I_1^Y = &\, Y \otimes w_0\,; \nn \\
I_2 =&\, \textstyle{\frac 12}f_{abc}\big(J^a \wedge J^b\big)\otimes w^c\,;\nn \\
I_2^W =&\, \textstyle{\frac 12}C_{abi}\big(J^a \wedge J^b\big)\otimes w^i\,;\nn \\
I_2^Y =&\, \big(Y\wedge J^a\big)\otimes V^a\,.
\end{align}
Here $f_{abc}$ is the totally antisymmetric structure tensor (assumed lowered with the $L_0$ Killing metric normalized to $\delta_{ab}$\,), and $C_{abi}$ is an invariant coupling between $W$ and the antisymmetric 
adjoint square. We evaluate the chain maps as follows:
\begin{align}
\partial_0  I_1 =&\, J^a  w^a = - f^{aa}{}_cw^c \equiv 0\,;\nn \\
 \partial_0  I_1^Y =&\, Yw_0 \equiv 0\,;\nn \\
\partial_1 I_2 =&\, f_{abc}\Big( J^a \otimes \big(J^b\otimes w^c\big)
 - \textstyle{\frac 12}f^{ab}{}_d J^d \otimes w^c\Big) =\nn \\
=&\, \big(f_{cab}f^{ab}{}_d  - \textstyle{\frac 12}f_{abc} f^{ab}{}_d\big) J^c \otimes w^d
 \equiv C \big(J^a \otimes w^a\big)\,;\nn \\
\partial_1 I_2^W \equiv &\, 0\,;\nn \\
\partial_1 I_2^Y =&\, Y \otimes (J^a w^a) - J^a \otimes ( Y w^a) \equiv 0\,.
\end{align}
Here the fact that $U \otimes U^*$ has zero $Y$-weight has been used. 
The evaluation of $\partial_1 I_2$ yields $I_1$ up to a multiple of the 
(nonzero) quadratic Casimir eigenvalue $C$ of $L_0$
in the adjoint representation\,. Further, because of 
the equivariance of  $\partial_1$\,, the image of $I_2^W$ must be an $L_0$ singlet,
which is excluded as there is, by 
construction\footnote{Explicitly, $\partial_1 I_2^W= C_{aj}J^a\otimes w^j$
where $C_{aj}:= C_{abi}D_j^{bi}$ in terms of the matrix elements of $J^b$ in the 
$\{w ^j \}$ basis.}, no invariant coupling between the adjoint
representation and a non-adjoint irreducible module $W$.  In summary, we have therefore
\[
\partial_0  I_1 =0= \partial_0  I_1^Y \,; \quad \partial_1 I_2 \,\propto\,  I_1\,; \quad
\partial_1 I_2^W =0\, \,\,\forall \,W\,, \quad \partial_1 I_2^Y=0\,,
\]
whence finally, as the singlet $Y \otimes w_0$ is the only nontrivial cycle within $B_1$, we have
\begin{align}
\mbox{dim}(\mbox{Ker} \partial_0) - \mbox{dim}(\mbox{Im} \partial_1) 
= 2-1 =\mbox{dim}({\mathbb C}I_1^Y) \equiv 1\,.
\end{align}
\mbox{} \hfill $\Box$\\

\section{Discussion}
\label{sec:Discussion}
The results reported in this paper address the  
structure of indecomposable representations of Lie superalgebras, specifically in establishing the existence of indecomposable doubles of Kac modules, for certain type I Lie superalgebras.

The methodology -- the computation of appropriate superalgebra (co)homologies, as elaborated in sections \ref{sec:KacModules}, \ref{sec:Relative} above -- is standard, once the ``double'' induction is set up (Lemma \ref{lem:DoubleSided} above\,).
In contrast to the semisimple Lie algebra case, it is evident that the theory is rather rich, even in finite dimensions. By the use of Shapiro's lemma (Lemma \ref{lem:Shapiro}), the computation is reduced to (co)homologies of (irreducible) modules of the even Lie subalgebra, which being reductive yields nontrivial results.


As stated above, the type I Lie superalgebras 
$sl(m/n)$\,, $osp(2/2n)$ belong to the basic classical class. We cite for example
\cite{Serganova1996KL,SuZhang2007character} and references therein, for results on the analysis of the Kac modules themselves, for general Lie superalgebras (the computation of irreducible factors in the composition series thereof, and thus evaluations of the characters of such modules). Here we have been at pains to provide direct insight into the existence of the indecomposable doubles,  irrespective of the explicit details of their construction, via specific (co)homology evaluations (in the atypical case, the Kac modules themselves are of course instances of other types of indecomposables). For definiteness, and given the central motivation coming from physics referred to in the introduction above, we confine the discussion here mainly to examples in $sl(2/1) \cong osp(2/2)$\,.  

Irreducible representations of $sl(2/1)$ are labelled in Kac-Dynkin notation by their highest weight $(a,b)$ with integer $a\ge 0$ the $sl(2)$ label of the highest weight (spin $j=\textstyle{\frac 12}a$\,),  and complex $b$\,. Typical modules have dimension $4(a+1)$ and comprise spin multiplets 
$j_y$\,, $(j\!\pm \!\textstyle{\frac 12})_{y-1}$\,, and $j_{y-2}$\,, labelled also by the additional Cartan generator $Y$ (`hypercharge') which commutes with $sl(2)$ (with eigenvalue $y=2b-a$ on the highest weight). Atypical modules (type 1) satisfy atypicality $b=0$ and in consequence have $y=-a$, dimension $2a+1$ and spin decomposition
$j_{-a}$ and $(j-\textstyle{\frac 12})_{-a-1}$\, while atypical modules (type 2)
satisfy atypicality $b=a+1$ and in consequence have $y=a+2$\,, dimension
$2a+3$ and spin decomposition $j_{a+2}$ and $(j+\textstyle{\frac 12})_{a+1}$\,.
For ease of notation, where no confusion arises we refer to these (irreducible) $sl(2/1)$\, modules by their dimension and leading hypercharge, $\underline{D}_y$\,: dimensions divisible by 4 being typical irreducible Kac modules (unless $y=-a$ or $y= a+2$), and odd dimensions being atypical irreducible of type 1 or type 2,
with the latter distinguished by an asterisk. Refer to table \ref{tab:EvenOddDecomp} for a listing of some low dimensional modules, and a summary of the notation used.

As discussed, construction of indecomposable module extensions requires
nonvanishing 1-(co)homology of $sl(2/1)$ valued in the appropriate tensor product of modules. Computation of 
$H^1$ valued in irreducible modules of $sl(2/1)$\, has been carried out in \cite{tanaka1995homology}, with results subsequently extended to all irreducibles of $sl(m/n)$ and $osp(2/2n)$ in \cite{su2007cohomology}\,.
In fact, for $sl(2/1)$\,, it is found \cite{tanaka1995homology,su2007cohomology} that the only irreducibles having nonvanishing $H^1$ are indeed the fundamental $\underline{3}_{-1}$\,, or its dual $\underline{3}^*_{2}$\,.  

These results are in fact 
informative about the structure of the simplest indecomposables: since $\underline{1}_{0}\otimes \underline{3}^*_{2} \cong \underline{3}^*_{2}$ has nonvanishing $H_1$, there exists
an extension, namely the atypical Kac module
$\underline{4}_{-1}\cong \underline{3}_{-1}\semiplus +\hskip5pt \underline{1}_{0}$\,. 
In a similar way we recover
the dual $\underline{4}_{2}$\, thanks to the corresponding trivial tensor product.
Likewise, the presence of $\underline{3}_{-1}$ or $\underline{3}^*_{2}$ in the resolution of the associated tensor products is responsible for example for
a shifted adjoint, atypical octet
$\underline{8}_{-1}\cong \underline{5}_{-2}\semiplus +\hskip5pt \underline{3}_{-1}$\,, and the dodecuplet 
$\underline{12}_{-3}\cong \underline{7}_{-3}\semiplus +\hskip5pt \underline{5}_{-2}$\,. 
By the same token $\underline{3}_{-1}\semiplus +\hskip5pt \underline{3}_{-1}$ is disallowed,
because $\underline{3}_{-1}\otimes \underline{3}^*_{2}
\cong \underline{8}_{1}+ \underline{1}_{0}$\,, with a similar conclusion
for  $\underline{7}_{-3}\semiplus +\hskip5pt \underline{3}_{-1}$ since
$\underline{7}_{-3}\otimes \underline{3}^*_{2}
\cong \underline{16}_{-1}+ \underline{5}_{-2}$\,. See table \ref{tab:TensorProdsSl21} for
a summary of these observations using tensor product resolutions in $sl(2/1)$ (see also \cite{marcu1980tensor})\,. Similar implications could also be drawn in more general cases, using the results of \cite{su2007cohomology}. 

In general, the existence of pairwise extensions $U \semiplus +\hskip3pt V$\,, 
$V \semiplus +\hskip3pt W$\, $\cdots$, is necessary (but not sufficient) for the construction of extended
indecomposable modules with structure $U\semiplus +\hskip3pt V \semiplus +\hskip3pt W \cdots$, in that the successive pairs are present as parts of appropriate factor spaces and invariant submodules. The directed acyclic graphs of such composition chains can also contain loops (for an example see below). The complexity of the problem of enumerating all indecomposables has been identified with that of classifying certain types of representations of free algebras \cite{germoni1997representations}. 

Finally we return to the case of indecomposable representations of $SU(2/1)$ with application to particle physics and the standard model. The fact that a single
family generation of quarks can be accommodated in an irreducible typical quartet representation $\underline{4}_{+2/3}$\,, with leptons in an irreducible $\underline{3}_{-1}$\,, (or
indecomposable (Kac module) $\underline{4}_{-1}$\,, if a right handed neutrino is present) was noted in the original papers
\cite{ne1979irreducible,fairlie1979higgs,dondi1979supersymmetric,neeman1980geometrical,thierry1982exterior}\,
(the hypercharges are in accord with standard assignments to give the correct fractional $+\textstyle{\frac 23},-\textstyle{\frac 13}$
quark electric charges, and $0,-1$ leptonic electric charges, through
the formula $Q=I_3 + \textstyle{\frac 12}Y$ in terms of weak isospin and hypercharge; see table \ref{tab:EvenOddDecomp}\,).
The indecomposable doubling 
$\underline{D}_y\semiplus +\hskip3pt \underline{D}_y$\,, for Kac modules
of $sl(m/n)$ and $osp(2/2n)$\,, the main result established in the present paper 
(Corollary \ref{cor:UniqueIndecDoubling} above) applies in particular to the doubling of the quartet representation, 
so that a single multiplet of $SU(2/1)$ can encompass both a first and a second family generation \cite{coquereaux1991elementary}, with consequent mixing of charged particle states (and violation of weak interaction 
universality through the coupling strength modification due to the Cabibbo mixing parameter $\sin \theta_C$ between $s$ and $d$ quarks\,), or of right-handed neutrino states (with consequent neutrino mass mixing and oscillation phenomena), respectively.
From the present perspective, the doubling necessitates a nonvanishing (co)homology in the tensor product of the appropriate quartet representations, namely $\underline{4}_{+2/3} \otimes \underline{4}_{+4/3}$\, (or for leptons,
$\underline{4}_{0} \otimes \underline{4}_{2}$\, ). This
indeed has dimension one, being a special case of the property established in this work. Thus, within $SU(2/1)$\,, the physically observed quantity $\sin \theta_C$ acquires an algebraic interpretation, as a parameter choice within 
the appropriate $H_1$ space. 

For the indecomposable doubling of quartet representations \cite{marcu1980representations}, the central 
case motivating the present work, it should be noted that the structure of the tensor product modules
$\underline{4}_{2/3} \otimes \underline{4}_{4/3}$\,
and $\underline{4}_{0} \otimes \underline{4}_{2}$\,, or more generally, of the family
$\underline{4}_{y} \otimes \underline{4}_{-y+2}$\, or $\underline{4}_{y} \otimes \underline{4}_{y'}$, has been studied by several authors in the course of their analyses of representations of
$sl(2/1)$\,
\cite{ScheunertNahmRittenberg1977irreps,marcu1980tensor, benamor1997reduction,ScheunertZhang1998}\,. In particular, in \cite{ScheunertZhang1998}\,, an explicit computation indeed yielded  (co)homology dimension 1, confirming our more general result. As noted above, in the particle physics context, this is consistent with the existence of the Cabibbo mixing angle $\theta_C$\,.
A more complete analysis of this family of tensor products and associated indecomposable modules, including cases where the composition chains form a loop, is presented in \cite{ThierryMiegJarvis2021indec}\,. 

A final point, of crucial significance for physics, is that the quartet double construction for $SU(2/1)$ has indeed been extended to triples 
$\underline{4}_y\semiplus +\hskip3pt \underline{4}_y \semiplus +\hskip3pt \underline{4}_y$\, (compare \cite{Su1992,ThierryMiegJarvis2021indec}) and beyond (for details see \cite{ThierryMiegJarvis2021indec}); a theory for general Kac modules is developed for $sl(n/1)$ and $sl(m/n)$ in \cite{germoni1997representations,
germoni1998indecomposable}.

In relation to the particle spectrum in the standard model, threefold quartet extensions are able to account for the experimental fact of 
three families of fundamental fermions (quarks and leptons, assigned to $SU(2/1)$ quartets), with the concomitant three generation charged quark electroweak mixing 
and modelling of the Cabibbo-Kobayashi-Maskawa matrix, as well as
neutrino mass mixing and oscillation phenomena
\cite{coquereaux1991elementary,ScheckHaussling1998triangular,haussling1998leptonic}\,.
 
\vfill
\noindent
\textbf{Acknowledgements:}\\
The authors wish to thank  Prof Ruibin Zhang for suggesting 
improvements to a preliminary version of the paper. We are grateful to Prof J\'{e}r\^{o}me Germoni for discussions and providing unpublished results. 
This research was supported in part by the Intramural Research Program of the National Library of Medicine, National Institute of Health.

\pagebreak
\begin{table}[tbp]
\begin{tabular}{|r|c|l|}
\hline \hline 
&&\\[-.15cm]
 $\underline{D}{}_y$ & $(a,b)$ & $\hskip1cm \sum\oplus j_y$ \\[-.15cm]
&&\\[-.15cm]
\hline
&&\\[-.15cm]
  $\underline{4}_y$ & $(0,b)$  &  $0_y+ \textstyle{\frac 12}_{y\!-\!1} + 0_{y\!-\!2}$\\[-.15cm]
&&\\[-.15cm]
\hline
&&\\[-.15cm]
  $\underline{8}_y$    &  $(1,b)$  &  $\textstyle{\frac 12}_{y}+ \big(0+1\big)_{y\!-\!1} + \textstyle{\frac 12}_{y\!-\!2}$\\[-.15cm]
&&\\[-.15cm]
\hline
&&\\[-.15cm]
  $\underline{8}_1$    &  $(1,1)$  &  $\textstyle{\frac 12}_{+1}+ \big(0+1\big)_{0} + \textstyle{\frac 12}_{-1}$\\[-.15cm]
&&\\[-.15cm]
\hline
&&\\[-.15cm]
  $\underline{12}_y$    &  $(2,b)$  &  $1_y+ (\textstyle{\frac 12}+\textstyle{\frac 32})_{y\!-\!1} + 1_{y\!-\!2}$    \\[-.15cm]
&&\\[-.15cm]\hline
&&\\[-.15cm]
  $\underline{4(a\!+\!1)}{}_y$  & $(a,b)$  &  $j_y +(j\pm\textstyle{\frac 12})_{y\!-\!1}
+   j_{y\!-\!2}$  \\[-.15cm]
&&\\[-.15cm]
\hline
\hline
\end{tabular}
\begin{tabular}{|l|c|l|l|}
\hline \hline 
&&&\\[-.15cm]
 $\underline{D}{}_y$ & $(a,b)$ & $\hskip.5cm \sum\oplus j_y$ & \\[-.15cm]
&&&\\[-.15cm]
\hline
&&&\\[-.15cm]
  $\underline{3}_{-1}$ & $(1,0)$  &  $\textstyle{\frac 12}_{-1}+  0_{-2}$ & 
$\scalebox{.7}{\raisebox{-0.1cm}{\young(~)}}$\\[-.15cm]
&&&\\[-.15cm]
\hline
&&&\\[-.15cm]
  $\underline{3}^*_2$    &  $(0,1)$  &  $0_{+2} +\textstyle{\frac 12}_{+1}$ & 
$\scalebox{.7}{\raisebox{-0.1cm}{$\overline{\young(~)}$}}$\\[-.15cm]
&&&\\[-.15cm]
\hline
&&&\\[-.15cm]
  $\underline{5}_{-2}$ & $(2,0)$  &  $1_{-2}+\textstyle{\frac 12}_{-3}$ & 
$\scalebox{.7}{\raisebox{-0.1cm}{${\young(~~)}$}}$\\[-.15cm]
&&&\\[-.15cm]
\hline
&&&\\[-.15cm]
  $\underline{5}^*_3$  & $(1,2)$  & $\textstyle{\frac 12}_{+3}+1_{+2}$ &
$\scalebox{.7}{\raisebox{-0.1cm}{$\overline{\young(~~)}$}}$\\[-.15cm]
&&&\\[-.15cm]
\hline
&&&\\[-.15cm]
  $\underline{7}_{-3}$    & $(3,0)$  & $\textstyle{\frac 32}_{-3}+1_{-4}$ &
$\scalebox{.7}{\raisebox{-0.1cm}{${\young(~~~)}$}}$\\[-.15cm]
&&&\\[-.15cm]
\hline
&&&\\[-.15cm]
  $\underline{7}^*_{4}$   & $(2,3)$  & $1_{+4}+\textstyle{\frac 32}_{+3}$ &
$\scalebox{.7}{\raisebox{-0.1cm}{$\overline{\young(~~~)}$}   }   $\\[-.15cm]
&&&\\[-.15cm]
\hline
\hline
\end{tabular}\mbox{}\\[.5cm]
\caption{\mbox{}\\
Left: $sl(2/1)$ typical irreducible representations with
highest weight Kac-Dynkin label $(a,b)$\,, giving alternative labelling
$\underline{D}_y$\, by dimension $D=4a\!+\!4$\,, with subscript the hypercharge $y=2b\!-\!a$ of the leading $sl(2)$ multiplet. The accompanying
decomposition into spin and hypercharge multiplets $j_{y}$
($j=\textstyle{\frac 12}a$\,) is given. The adjoint is $\underline{8}_1$\,. \mbox{}
\\
Right: $sl(2/1)$ representations and spin hypercharge decompositions for atypical
irreducibles $(a,b)$\, of types 1 ( unstarred) and 2 ( starred${}^*$\,).
See text for details of the notation. 
Associated Young diagrams are also provided {\protect\cite{dondi1981diagram}}\,.
\protect\label{tab:EvenOddDecomp} 
}
\end{table}

\begin{table}
\begin{tabular}{|rcl|rcl|}
\hline \hline 
  &&&&&\\[-.15cm]
 $\underline{D}\otimes \underline{D}'$ & $\cong$  & $ \sum\oplus \underline{D}''$ & 
&  &\\[-.15cm]
  &&&&&\\[-.15cm]
\hline
\hline
  &&&&&\\[-.15cm]
$\underline{3}\otimes \underline{3}$  & $\cong$ & $\underline{5}+ \underline{4}$ & 
$\scalebox{.7}{\raisebox{-0.1cm}{\young(~)}}\otimes
\scalebox{.7}{\raisebox{-0.1cm}{\young(~)}}$& $\cong$&
$\scalebox{.7}{\raisebox{-0.1cm}{\young(~~)}} +
\scalebox{.7}{\raisebox{-0.3cm}{\young(~,~)}}$\\[-.15cm]
  &&&&&\\[-.15cm]
\hline
  &&&&&\\[-.15cm]
$\underline{3}\otimes \underline{3}^*$  & $\cong$ & $\underline{8}+ \underline{1}$ & 
$\scalebox{.7}{\raisebox{-0.1cm}{${\young(~)}$}}\otimes \scalebox{.7}{\raisebox{-0.1cm}{$\overline{\young(~)}$}}$& $\cong$&
$\scalebox{.7}{\raisebox{-0.1cm}{$\overline{\young(~)}\young(~)$}} +
\circ$ \\[-.15cm]
  &&&&&\\[-.15cm]
\hline
  &&&&&\\[-.15cm]
$\underline{5}\otimes \underline{3}^*$  & $\cong$ &$\underline{12} 
+ \underline{3}$ & $\scalebox{.7}{\raisebox{-0.1cm}{$\young(~~)$}}\otimes
\scalebox{.7}{\raisebox{-0.1cm}{$\overline{\young(~)}$}}$ & $\cong$ &
$\scalebox{.7}{\raisebox{-0.1cm}{$\overline{\young(~)}\young(~~)$}} +
\scalebox{.7}{\raisebox{-0.1cm}{$\young(~)$}}$\\[-.15cm]
  &&&&&\\[-.15cm]
\hline
  &&&&&\\[-.15cm]
$\underline{7}\otimes \underline{3}^*$  & $\cong$ &$\underline{16} 
+ \underline{5}$& $\scalebox{.7}{\raisebox{-0.1cm}{$\young(~~~)$}}\otimes
\scalebox{.7}{\raisebox{-0.1cm}{$\overline{\young(~)}$}}$& $\cong$&
$\scalebox{.7}{\raisebox{-0.1cm}{$\overline{\young(~)}\young(~~~)$}} +
\scalebox{.7}{\raisebox{-0.1cm}{${\young(~~)}$}}$\hskip.3cm\\[-.15cm]
  &&&&&\\[-.15cm]
\hline
\hline
\end{tabular}
\mbox{}\\[.5cm]
\caption{Some $sl(2/1)$ tensor product resolutions with
associated Young diagrams \cite{dondi1981diagram}\,.
See text for details.
\label{tab:TensorProdsSl21}}
\end{table}
\mbox{}
\vfill
\newpage

\begin{thebibliography}{10}
\bibliographystyle{plain}

\bibitem{Kac1977classification}
V.~G. Kac.
\newblock Lie superalgebras.
\newblock {\em Advances in Math.}, 26(1):8--96, 1977.

\bibitem{Kac1977characters}
V.~G. Kac.
\newblock Characters of typical representations of classical {L}ie
  superalgebras.
\newblock {\em Comm. Algebra}, 5(8):889--897, 1977.

\bibitem{scheunert2006theory}
Manfred Scheunert.
\newblock {\em The theory of Lie superalgebras: an introduction}, volume 716.
\newblock Springer, 2006.

\bibitem{frappat2000dictionary}
Luc Frappat, Antonino Sciarrino, and Paul Sorba.
\newblock {\em Dictionary on Lie algebras and superalgebras}, volume~10.
\newblock Academic Press San Diego, CA, 2000.

\bibitem{cornwell1992group}
John~F Cornwell.
\newblock {\em Group theory in physics. Volume {III} {S}upersymmetries and
  infinite-dimensional algebra ({T}echniques of {P}hysics)}.
\newblock Academic Press, 1992.

\bibitem{ScheunertNahmRittenberg1977irreps}
M.~Scheunert, W.~Nahm, and V.~Rittenberg.
\newblock Irreducible representations of the {${\rm osp}(2,1)$} and {${\rm
  spl}(2,1)$} graded {L}ie algebras.
\newblock {\em Journal of Mathematical Physics}, 18(1):155--162, 1977.

\bibitem{marcu1980tensor}
Mihael Marcu.
\newblock The tensor product of two irreducible representations of the
  $spl(2,1)$ superalgebra.
\newblock {\em Journal of Mathematical Physics}, 21(6):1284--1292, 1980.

\bibitem{marcu1980representations}
Mihael Marcu.
\newblock The representations of $spl(2,1)$ -- an example of representations of
  basic superalgebras.
\newblock {\em Journal of Mathematical Physics}, 21(6):1277--1283, 1980.

\bibitem{Su1992}
Yucai Su.
\newblock Classification of finite dimensional modules of the {L}ie
  superalgebra sl(2/1).
\newblock {\em Communications in Algebra}, 20(11):3259--3277, 1992.

\bibitem{gotz2007representation}
Gerhard G{\"o}tz, Thomas Quella, and Volker Schomerus.
\newblock Representation theory of $sl(2|1)$.
\newblock {\em Journal of Algebra}, 312(2):829--848, 2007.

\bibitem{germoni1997representations}
J{\'e}r{\^o}me Germoni.
\newblock Repr{\'e}sentations ind{\'e}composables des superalgebres de {L}ie
  sp{\'e}ciales lin{\'e}aires $sl(m/1)$.
\newblock {\em Comptes Rendus de l'Acad{\'e}mie des Sciences-Series
  I-Mathematics}, 324(11):1221--1226, 1997.

\bibitem{germoni1998indecomposable}
J{\'e}r{\^o}me Germoni.
\newblock Indecomposable representations of special linear Lie superalgebras.
\newblock {\em Journal of Algebra}, 209(2):367--401, 1998.

\bibitem{germoni2000indecomposable}
J{\'e}r{\^o}me Germoni.
\newblock Indecomposable representations of $osp(3,2)$\,, $ {D} (2, 1;\alpha)$
  and ${G}(3)$.
\newblock {\em Boletin de la Academia Nacional de Ciencias}, 65:147, 2000.

\bibitem{Leites1975cohomology}
D.~A. Le\u{i}tes.
\newblock Cohomology of {L}ie superalgebras.
\newblock {\em Funkcional. Anal. i Prilo\v{z}en.}, 9(4):75--76, 1975.

\bibitem{Serganova1996KL}
Vera Serganova.
\newblock Kazhdan-{L}usztig polynomials and character formula for the {L}ie
  superalgebra {${\mathfrak g}{\mathfrak l}(m|n)$}.
\newblock {\em Selecta Math. (N.S.)}, 2(4):607--651, 1996.

\bibitem{SuZhang2007character}
Yucai Su and R.~B. Zhang.
\newblock Character and dimension formulae for general linear superalgebra.
\newblock {\em Adv. Math.}, 211(1):1--33, 2007.

\bibitem{ne1979irreducible}
Yuval Ne'eman.
\newblock Irreducible gauge theory of a consolidated {S}alam-{W}einberg model.
\newblock {\em Physics Letters B}, 81(2):190--194, 1979.

\bibitem{fairlie1979higgs}
D~B Fairlie.
\newblock Higgs fields and the determination of the {W}einberg angle.
\newblock {\em Physics Letters B}, 82(1):97--100, 1979.

\bibitem{dondi1979supersymmetric}
P~H Dondi and Peter~D Jarvis.
\newblock A supersymmetric {W}einberg-{S}alam model.
\newblock {\em Physics Letters B}, 84(1):75--78, 1979.

\bibitem{neeman1980geometrical}
Yuval Ne'eman and Jean Thierry-Mieg.
\newblock Geometrical gauge theory of ghost and {G}oldstone fields and of ghost
  symmetries.
\newblock {\em Proceedings of the National Academy of Sciences},
  77(2):720--723, 1980.

\bibitem{thierry1982exterior}
Jean Thierry-Mieg and Yuval Ne'eman.
\newblock Exterior gauging of an internal supersymmetry and {SU(2/1)} quantum
  asthenodynamics.
\newblock {\em Proceedings of the National Academy of Sciences},
  79(22):7068--7072, 1982.

\bibitem{ThierryMieg2021newchi}
Jean Thierry-Mieg.
\newblock Chirality, a new key for the definition of the connection and
  curvature of a {L}ie-{K}ac superalgebra.
\newblock {\em Journal of High Energy Physics}, 2021(1):111, 2021.

\bibitem{ThierryMiegJarvis2021newchi}
Jean Thierry-Mieg and Peter Jarvis.
\newblock ${S}{U}(2/1)$ superchiral self-duality: a new quantum, algebraic and
  geometric paradigm to describe the electroweak interactions.
\newblock {\em Journal of High Energy Physics}, 2021(4):1, 2021.

\bibitem{coquereaux1991elementary}
Robert Coquereaux.
\newblock Elementary fermions and $su(2|1)$ representations.
\newblock {\em Physics Letters B}, 261(4):449--458, 1991.

\bibitem{ScheckHaussling1998triangular}
Rainer H\"aussling and Florian Scheck.
\newblock Triangular mass matrices of quarks and
  {C}abibbo-{K}obayashi-{M}askawa mixing.
\newblock {\em Phys. Rev. D}, 57:6656--6662, Jun 1998.

\bibitem{haussling1998leptonic}
R~Haussling, M~Paschke, and F~Scheck.
\newblock Leptonic generation mixing, noncommutative geometry and solar
  neutrino fluxes.
\newblock {\em Physics Letters B}, 417(3-4):312--319, 1998.

\bibitem{ThierryMiegJarvis2021indec}
Jean Thierry-Mieg, Peter Jarvis and J{\'e}r{\^o}me Germoni.
\newblock Explicit construction of the finite dimensional indecomposable
  representations of the simple {L}ie-{K}ac ${S}{U}(2/1)$ superalgebra and
  their low level non diagonal super {C}asimir operators.
\newblock {P}reprint, arXiv:2207.06545;
\newblock Jean Thierry-Mieg, Peter Jarvis and J{\'e}r{\^o}me Germoni.
\newblock Construction of indecomposable $N$-replications of Kac modules of type 1 Lie superalgebras $sl(m/n), osp(2/2n)$.
\newblock {P}reprint, arXiv:2207.06538.


\bibitem{BenAmorPinczon1991}
Hedi Ben~Amor and Georges Pinczon.
\newblock Extensions of representations of {L}ie superalgebras.
\newblock {\em Journal of Mathematical Physics}, 32(3):621--629, 1991.

\bibitem{ChevalleyEilenberg1948}
Claude Chevalley and Samuel Eilenberg.
\newblock Cohomology theory of {L}ie groups and {L}ie algebras.
\newblock {\em Trans. Amer. Math. Soc.}, 63:85--124, 1948.

\bibitem{ScheunertZhang1998}
M.~Scheunert and R.~B. Zhang.
\newblock Cohomology of {L}ie superalgebras and their generalizations.
\newblock {\em Journal of Mathematical Physics}, 39(9):5024--5061, 1998.

\bibitem{tanaka1995homology}
Junko Tanaka.
\newblock Homology and cohomology of {L}ie superalgebra $sl(2/1)$ with
  coefficients in the spaces of finite-dimensional irreducible representations.
\newblock {\em Journal of Mathematics of Kyoto University}, 35(4):733--756,
  1995.

\bibitem{su2007cohomology}
Yucai Su and R~B Zhang.
\newblock Cohomology of {L}ie superalgebras $\mathfrak{sl}(m|n)$ and
  $\mathfrak{osp}(2|2n)$.
\newblock {\em Proceedings of the London Mathematical Society}, 94(1):91--136,
  2007.

\bibitem{benamor1997reduction}
Hedi Benamor.
\newblock Reduction of the adjoint representation of $sl(2,1)$: Generators of
  primitive ideals.
\newblock {\em Communications in algebra}, 25(3):715--736, 1997.

\bibitem{dondi1981diagram}
P~H Dondi and P~D Jarvis.
\newblock Diagram and superfield techniques in the classical superalgebras.
\newblock {\em Journal of Physics A: Mathematical and General}, 14(3):547,
  1981.

\end{thebibliography}

\vfill
\begin{appendix}
\section{Double induction equivalence} 
\label{sec:DoubleSided}
\noindent

\noindent
\textbf{Proof of Lemma \ref{lem:DoubleSided}:}\\
Let us write 
$V$ for the highest weight Kac module 
\smash{$\mbox{Ind}\! \left.\right|_{L_{+}}^L (U)=\overline{V_+(U)}{}$},
and $W$ for the $L$-module \smash{$\mbox{Ind}\! \left.\right|_{L_{\zbar}}^L \big(U \otimes U^*\big)$} induced from the $L_{\zbar}$-module 
$U \otimes U^*$\,. The map $\pi: V\otimes V^* \rightarrow W$ confers a twisted coproduct action on $W$\,, with generators defined by
\begin{align}
\widehat{J} := \pi \circ \Delta J \circ \pi^{-1}\,,
\quad \widehat{Q} := \pi \circ \Delta Q\circ \pi^{-1}\,,
\quad \widehat{\Qbar} := \pi \circ \Delta \Qbar \circ \pi^{-1}\,,
\end{align}
which we wish to show is equivalent to the direct (left) action on the induced module $W$\,. On states $w = m \om \otimes (u\otimes u^*)$ we have for example
for even generators
\begin{align}
\Delta J \pi^{-1} (w) =&\, (Jm \otimes_+ u)\otimes (\om \otimes_- u^*)
+ (m \otimes_+ u)\otimes (J\om \otimes_- u^*) \nn \\
=&\, ([J,m] \otimes_+ u)\otimes (\om \otimes_- u^*)
+ (m \otimes_+ Ju)\otimes (\om \otimes_- u^*)+\nn \\
&\, +
(m \otimes_+ u)\otimes ([J,\om] \otimes_- u^*)
+ (m \otimes_+ u)\otimes (\om \otimes_- Ju^*)\,,\nn \\
\mbox{so}\quad
\pi \Delta J \pi^{-1} w = &\,[J,m\om]\otimes (u \otimes u^*)+
m\om\otimes \Delta J \cdot(u \otimes u^*) \nn \\
= &\,  J m\om\otimes (u \otimes u^*)\,.
\end{align}
Thus 
we have $\widehat{J}w \equiv J w$\,. On the other hand\footnote{Assigning $m\om$ to be of degree $(p,q)$ as a monomial in odd generators\,.}
\begin{align}
\label{eq:QQbarAction}
\widehat{Q} w =&\, Qm \om \otimes (u\otimes u^*) + 
\sum m\om_Q \otimes (u \otimes  J_Q u^*)\,, \nonumber\\
\widehat{\Qbar} w =&\, 
\sum m_{\Qbar} \om \otimes  J_{\Qbar}u \otimes u^*
+ (-1)^p m \Qbar \om\otimes w\,, 
\end{align}
where, for example, (c.f. the evaluation of $\widehat{J} w$\,),
\begin{align*}
\Qbar m \otimes_+ u^* = &\, [\Qbar, m]\otimes_+ u^*
= \sum m_{\Qbar} J_{\Qbar} \otimes_+ u^*\equiv  
\sum m_{\Qbar}\otimes_+  J_{\Qbar} u^* \,,
\end{align*}
with the $ [\Qbar, m]$ bracket having been re-ordered as a sum of 
monomials with even generators moved to the right.
These additional terms 
thus constitute the difference between the twisted coproduct and standard actions on $W$ under the odd generators, whereas the even actions coincide (
$V\otimes V^*$ and $W$ are equivalent as $L_{\zbar}$-modules). 

From (\ref{eq:QQbarAction}) the $\widehat{Q}$-\,, $\widehat{\Qbar}$-actions shift the monomial gradings
from states in sector $(p,q)$ to states in sectors $(p+1, q)$\, $(p, q+1)$\,, respectively,
but with lower degree contributions $(p, q-1)$ and $(p-1, q)$\,.
In order to establish the equivalence with the standard actions it is necessary to find an additional straightening map $\varphi:W \rightarrow W$\,, which acts within sectors $(p,q) + (p-1,q-1) + (p-2,q-2) +\cdots$ to remove trailing terms. 

We now show this construction explicitly for the case of $sl(2/1)$\,. We take an $sl(2)$-covariant basis
 for $L=sl(2/1)$ with even generators 
$J_i$, $Y$\,, $i=1,2,3$\,  and odd generators (lowering and raising operators)
$Q_a$, $\Qbar{}^b$, $a,b = 1,2$\,. The (anti)commutation relations are
\footnote{Here $\mbox{\boldmath{$\sigma$}}=(\sigma_1, \sigma_2,\sigma_3)$
are the standard $2\times 2$ Pauli matrices. We  use the 
antisymmetric symbols  $\varepsilon_{12}=1 = -\varepsilon^{12}$\,,
$\varepsilon^{ab} \varepsilon_{bc}:= \delta^a{}_c$\,, and 
$\varepsilon_{ijk}$ with $\varepsilon_{123}=1$\,. }
\begin{align}
{[}J_i, J_j{]} = &\, \varepsilon_{ijk} J_k\,, \quad [Y, J_i] =0\,, \nn \\
{[}J_i, Q_a{]} = &\, -\textstyle{\frac 12}Q_b\big(\sigma_i)^b{}_a\,,\quad
{[}J_i, \Qbar{}^a{]} = +\textstyle{\frac 12}\big(\sigma_i)^a{}_b\Qbar{}^b\,,\nn \\
[Y, Q]=&\,-Q\,, \quad [Y, \Qbar{}^a]=+\Qbar{}^a\,, \nn \\
{\{}\Qbar{}^a,Q_b\}=&\, Y\delta^a{}_b - 
\big(\mbox{\boldmath{$\sigma$}}\cdot{\mathbf J}\big){}^a{}_b\,,\nn \\
\mbox{with}\qquad {\{}Q_a, Q_b{\}}=&\, 0 = {\{}\Qbar{}^a,\Qbar{}^b{\}} \,.
\end{align}

There are 9 $(p,q)$ sectors, $0\le p,q \le 2$, of the form $m \om \otimes w$, with 10 spin projections, including both spin-1 and spin-0 parts within $(1,1) \cong (1,1)_1+(1,1)_0$\,.  Mappings between sectors are determined by $sl(2)$ coupling coefficients, up to reduced matrix elements. In tensor notation, this is taken into account by using appropriately symmetrized basis states. Let the $L_{\zbar} = gl(2)$-module $U$ have hypercharge $Y=y$\,, and spin $j$ with generators ${\mathbf j}$\,.
To simplify the notation we write $\otimes (u\otimes u^*)$ as $|\mbox{\hskip.2cm}\rangle $\,.
%
In spin-zero combinations the $sl(2)$ spin generators do not contribute,
as $Tr(\mbox{\boldmath{$\sigma$}})=0$\,, whereas 
in spin-1 projections the full $(\mbox{\boldmath{$\sigma$}}\!\cdot\!{\mathbf J}){}^a{}_b$ operator is present in matrix elements. Consider for example
\begin{align}
(\widehat{Q} \cdot \Qbar|\mbox{\hskip.2cm}\rangle)_0\equiv &\,
\widehat{Q}_a \cdot (\Qbar{}^a|\mbox{\hskip.2cm}\rangle)
 =({Q}\Qbar)\mh + 2(-y)\mh \,;\nn \\
\widehat{Q}_a\cdot (\Qbar\Qbar)\mh =&\,
{Q}_a((\varepsilon \Qbar{})_b\Qbar{}^b)\mh + 
(\textstyle{\frac 12}(\varepsilon \Qbar)_a + 
2 (\varepsilon \Qbar)_b (\mbox{\boldmath{$\sigma$}}\!\cdot\!{\mathbf j})^b{}_a)\mh\,,
\end{align}
where the scalars $(Q \Qbar)$, $(\Qbar \Qbar)$ are traces with respect to the $sl(2)$ invariants $\delta$ and $\varepsilon$, respectively.
After re-ordering, the even generators $Y$, ${\mathbf J}$ have been replaced by their representatives $(-y), {\mathbf j}$\,, acting on the second tensor factor of the state $u \otimes u^*$\,.

The appropriate matrix elements of the straightening map $\varphi$ are
therefore
\begin{align}
\varphi(\Qbar \mh)= &\,\Qbar \mh\,, \nn \\ 
\varphi((\Qbar\Qbar) \mh)= &\,
(\Qbar\Qbar) \mh\,, \nn \\
\varphi((Q\Qbar )\mh) =&\, (Q\Qbar )\mh +2y \mh\,,\nn \\
\varphi(Q(\Qbar \Qbar)\mh) = &\,
Q(\Qbar \Qbar)\mh - (\textstyle{\frac 12}(\varepsilon \Qbar) + 
2 (\varepsilon \Qbar)  \mbox{\boldmath{$\sigma$}}\!\cdot\! {\mathbf j})\mh\,.
\end{align}
The construction of $\varphi$ can be completed by defining it in a similar way on the remaining states which mix, namely in sector $(QQ)(\Qbar\Qbar)\mh$ with $(Q\Qbar) \mh$
and $\mh$\,, and in sector $Q(Q \Qbar)\mh$ with $Q\mh$\,, such that we have 
$\varphi \widehat{Q}\varphi^{-1} = Q$\, as required\,.

Turning to the equivalence of  $\Qbar$ and \smash{$\widehat{\Qbar}$} on $W$, a map $\overline{\varphi}$ analogous  to $\varphi$ can be constructed\footnote{Due to the ordering choice, construction for $\Qbar$ in the $m \om \otimes w$ basis would be more involved.}, but with a change of monomial ordering
to $\om m \otimes w$. However, the switch mapping is simply a nonsingular basis change $\sigma: W \rightarrow W$\,,
so that \smash{$\overline{\varphi} \widehat{\Qbar}\overline{\varphi}^{-1}$} and $\Qbar$ are brought into coincidence in the original basis, by a further conjugation with respect to $\sigma$\,.

For the general case of Lie superalgebras $sl(m/n)$ or $osp(2/2n)$\,, the construction proceeds in an analogous way, although the details are more technical. For example
even in $sl(3/1)$ there are now $16$ sectors $(p,q)$\,, with projections on to two irreducibles  if $p,q \ne 0$. However, it is clear that the 
required automorphism of $W$ (the straightening map $\varphi$), can be covariantly constructed incrementally in the same fashion, to be of lower triangular form with respect to the polynomial degree $p+q$\,, and with unit diagonal in a suitable basis.
\mbox{}\hskip1cm
\,\, \hfill $\Box$

\end{appendix}
\vfill
\end{document}